# ON NEW TYPES OF WEAKLY NANO CONTINUITY

QAYS HATEM IMRAN[1]

***Abstract.*** *In this paper, we shall use the concepts of* $N\alpha$*-open and* $NS\alpha$*-open sets to define some new types of weakly nano continuity such as;* $N\alpha$*-continuous,* $N\alpha$**-continuous,* $N\alpha$***-continuous,* $NS\alpha$*-continuous,* $NS\alpha$**-continuous and* $NS\alpha$***-continuous maps. Also, we shall explain the relationships between these types of weakly nano continuity and the concepts of nano continuity. Moreover, we shall prove some theorems, properties, remarks and give counter examples about these new concepts of weakly nano continuity.*

## 1. INTRODUCTION

M.L. Thivagar and C. Richard [1] present nano topological space (or simply $\mathrm{n.t.s.}$) on a subset $\mathcal{M}$ of a universe which is defined with respect to lower and upper approximations of $\mathcal{M}$. He studied about the weak forms of nano open sets. Q.H. Imran [3] presented the concept of NSα-open sets in nano topological spaces. The objective of this paper is to present new types of weakly nano continuity such as; Nα-continuous, Nα*-continuous, Nα**-continuous, NSα-continuous, NSα*-continuous and NSα**-continuous maps. Also, we shall explain the relationships between these types of weakly nano continuity and the concepts of nano continuity. Moreover, we shall prove some theorems, properties, remarks and give counter examples about these new concepts of weakly nano continuity.

## 2. PRELIMINARIES

Throughout this paper, $(\mathcal{U},\tau_{\mathcal{R}}(\mathcal{M})),(\mathcal{V},\sigma_{\mathcal{R}}(\mathcal{N}))$ and $(\mathcal{W},\rho_{\mathcal{R}}(\mathcal{O}))$ (or simply $\mathcal{U},\mathcal{V}$ and $\mathcal{W}$) constantly mean $\mathrm{n.t.s.}$ on which no separation axioms are normal unless for the most part determined. For a set $\mathcal{D}$ in a $\mathrm{n.t.s.}$ $(\mathcal{U},\tau_{\mathcal{R}}(\mathcal{M}))$, $\mathrm{Ncl}(\mathcal{D})$, $\mathrm{Nint}(\mathcal{D})$ and $\mathcal{D}^c=\mathcal{U}-\mathcal{D}$ denote the nano closure of $\mathcal{D}$, the nano interior of $\mathcal{D}$ and the nano complement of $\mathcal{D}$ respectively.

[1] Muthanna University, College of Education for Pure Science, Department of Mathematics, Iraq.
E-mail: qays.imran@mu.edu.iq.

**Definition 2.1.** A subset $\mathcal{D}$ of a $\text{n.t.s.}$ $(\mathcal{U}, \tau_{\mathcal{R}}(\mathcal{M}))$ is said to be:

i. A nano α-open set (in short Nα-open set) [1] if $\mathcal{D} \subseteq \text{Nint}(\text{Ncl}(\text{Nint}(\mathcal{D})))$. The family of all Nα-open sets of $\mathcal{U}$ is denoted by $\tau_{\mathcal{R}}\alpha(\mathcal{M})$.
ii. A nano semi-α-open set (in short NSα-open set) [3] if there exists a Nα-open set $\mathcal{P}$ in $\mathcal{U}$ such that $\mathcal{P} \subseteq \mathcal{D} \subseteq \text{Ncl}(\mathcal{P})$ or equivalently if $\mathcal{D} \subseteq \text{Ncl}(\text{Nint}(\text{Ncl}(\text{Nint}(\mathcal{D}))))$. The family of all NSα-open sets of $\mathcal{U}$ is denoted by $\tau_{\mathcal{R}}S\alpha(\mathcal{M})$. The complement of NSα-open set is called a nano semi-α-closed set (in short NSα-closed set).

**Example 2.2.** Let $\mathcal{U} = \{r_1, r_2, r_3, r_4\}$ with $\mathcal{U}/\mathcal{R} = \{\{r_1\}, \{r_3\}, \{r_2, r_4\}\}$ and $\mathcal{M} = \{r_1, r_2\}$. Then $\tau_{\mathcal{R}}(\mathcal{M}) = \{\phi, \{r_1\}, \{r_2, r_4\}, \{r_1, r_2, r_4\}, \mathcal{U}\}$ is a $\text{n.t.s.}$.
The family of all Nα-open sets of $\mathcal{U}$ is: $\tau_{\mathcal{R}}\alpha(\mathcal{M}) = \{\phi, \{r_1\}, \{r_2, r_4\}, \{r_1, r_2, r_4\}, \mathcal{U}\}$. The family of all NSα-open sets of $\mathcal{U}$ is: $\tau_{\mathcal{R}}S\alpha(\mathcal{M}) = \tau_{\mathcal{R}}\alpha(\mathcal{M}) \cup \{\{r_1, r_3\}, \{r_2, r_3, r_4\}\}$.

**Remark 2.3 [3].** In a $\text{n.t.s.}$ $(\mathcal{U}, \tau_{\mathcal{R}}(\mathcal{M}))$, then the following statements hold and the opposite of each statement is not true:

i. Every N-open set is a Nα-open and NSα-open.
ii. Every Nα-open set is a NSα-open.

**Definition 2.4.** Let $h: (\mathcal{U}, \tau_{\mathcal{R}}(\mathcal{M})) \longrightarrow (\mathcal{V}, \sigma_{\mathcal{R}}(\mathcal{N}))$ be a map, then $h$ is said to be:

i. Nano continuous (in short N-continuous) [2] iff for each $\mathcal{D}$ N-open set in $\mathcal{V}$, then $h^{-1}(\mathcal{D})$ is a N-open set in $\mathcal{U}$.
ii. Nano α-continuous (in short Nα-continuous) [4] iff for each $\mathcal{D}$ N-open set in $\mathcal{V}$, then $h^{-1}(\mathcal{D})$ is a Nα-open set in $\mathcal{U}$.

**Theorem 2.5 [2].** A map $h: (\mathcal{U}, \tau_{\mathcal{R}}(\mathcal{M})) \longrightarrow (\mathcal{V}, \sigma_{\mathcal{R}}(\mathcal{N}))$ is N-continuous iff $h^{-1}(\text{Nint}(\mathcal{D})) \subseteq \text{Nint}(h^{-1}(\mathcal{D}))$ for every $\mathcal{D} \subseteq \mathcal{V}$.

**Definition 2.6 [2].** Let $h: (\mathcal{U}, \tau_{\mathcal{R}}(\mathcal{M})) \longrightarrow (\mathcal{V}, \sigma_{\mathcal{R}}(\mathcal{N}))$ be a map, then $h$ is said to be nano open (in short N-open) iff for each $\mathcal{D}$ N-open set in $\mathcal{U}$, then $h(\mathcal{D})$ is a N-open set in $\mathcal{V}$.

## 3. WEAKLY NANO CONTINUOUS MAPS

**Definition 3.1.** Let $h: (\mathcal{U}, \tau_{\mathcal{R}}(\mathcal{M})) \longrightarrow (\mathcal{V}, \sigma_{\mathcal{R}}(\mathcal{N}))$ be a map, then $h$ is said to be:

i. Nano α*-continuous (in short Nα*-continuous) iff for each $\mathcal{D}$ Nα-open set in $\mathcal{V}$, then $h^{-1}(\mathcal{D})$ is a Nα-open set in $\mathcal{U}$.
ii. Nano α**-continuous (in short Nα**-continuous) iff for each $\mathcal{D}$ Nα-open set in $\mathcal{V}$, then $h^{-1}(\mathcal{D})$ is a N-open set in $\mathcal{U}$.

**Definition 3.2.** Let $h: (\mathcal{U}, \tau_{\mathcal{R}}(\mathcal{M})) \longrightarrow (\mathcal{V}, \sigma_{\mathcal{R}}(\mathcal{N}))$ be a map, then $h$ is said to be:

i. Nano semi-α-continuous (in short NSα-continuous) iff for each $\mathcal{D}$ N-open set in $\mathcal{V}$, then $h^{-1}(\mathcal{D})$ is a NSα-open set in $\mathcal{U}$.
ii. Nano semi-α*-continuous (in short NSα*-continuous) iff for each $\mathcal{D}$ NSα-open set in $\mathcal{V}$, then $h^{-1}(\mathcal{D})$ is a NSα-open set in $\mathcal{U}$.
iii. Nano semi-α**-continuous (in short NSα**-continuous) iff for each $\mathcal{D}$ NSα-open set in $\mathcal{V}$, then $h^{-1}(\mathcal{D})$ is a N-open set in $\mathcal{U}$.

**Theorem 3.3.** Let $h: (\mathcal{U}, \tau_{\mathcal{R}}(\mathcal{M})) \longrightarrow (\mathcal{V}, \sigma_{\mathcal{R}}(\mathcal{N}))$ be a map. Then the following statements are equivalent:

i. $h$ is a NSα-continuous.
ii. The inverse image of each N-closed set in $\mathcal{V}$ is NSα-closed set in $\mathcal{U}$.
iii. $h(\text{Nint}(\text{Ncl}(\text{Nint}(\text{Ncl}(\mathcal{C}))))) \subseteq \text{Ncl}(h(\mathcal{C}))$, for each $\mathcal{C} \in \mathcal{U}$.
iv. $\text{Nint}(\text{Ncl}(\text{Nint}(\text{Ncl}(h^{-1}(\mathcal{D}))))) \subseteq h^{-1}(\text{Ncl}(\mathcal{D}))$, for each $\mathcal{D} \in \mathcal{V}$.

*Proof:*
(i) ⟹ (ii). Let $\mathcal{D}$ be N-closed set in $\mathcal{V}$. This implies that $\mathcal{V} - \mathcal{D}$ is a N-open set.
Hence $h^{-1}(\mathcal{V} - \mathcal{D})$ is a NSα-open set in $\mathcal{U}$. i.e., $\mathcal{U} - h^{-1}(\mathcal{D})$ is a NSα-open set in $\mathcal{U}$.
Thus $h^{-1}(\mathcal{D})$ is a NSα-closed set in $\mathcal{U}$.
(ii) ⟹ (iii). Let $\mathcal{C} \in \mathcal{U}$, then $\text{Ncl}(h(\mathcal{C}))$ is a N-closed set in $\mathcal{V}$. So that $h^{-1}(\text{Ncl}(h(\mathcal{C})))$ is NSα-closed set in $\mathcal{U}$. Thus we have
$h^{-1}(\text{Ncl}(h(\mathcal{C}))) \supseteq \text{Nint}(\text{Ncl}(\text{Nint}(\text{Ncl}(h^{-1}(\text{Ncl}(h(\mathcal{C}))))))) \supseteq \text{Nint}(\text{Ncl}(\text{Nint}(\text{Ncl}(\mathcal{C}))))$.
Or $\text{Ncl}(h(\mathcal{C}))) \supseteq h(\text{Nint}(\text{Ncl}(\text{Nint}(\text{Ncl}(\mathcal{C})))))$.
(iii) ⟹ (iv). Since $\in \mathcal{V}$ , $h^{-1}(\mathcal{D}) \in \mathcal{U}$ so by hypothesis we have
$\text{Nint}(\text{Ncl}(\text{Nint}(\text{Ncl}(h^{-1}(\mathcal{D}))))) \subseteq \text{Ncl}(h(h^{-1}(\mathcal{D}))) \subseteq \text{Ncl}(\mathcal{D})$, that is
$\text{Nint}(\text{Ncl}(\text{Nint}(\text{Ncl}(h^{-1}(\mathcal{D}))))) \subseteq h^{-1}(\text{Ncl}(\mathcal{D}))$.
(iv) ⟹ (i). Let $\mathcal{C}$ be a N-open subset of $\mathcal{V}$. Let $\mathcal{D} = \mathcal{V} - \mathcal{C}$ and $\mathcal{C} = h^{-1}(\mathcal{D})$ by (iii) we have
$\text{Nint}(\text{Ncl}(\text{Nint}(\text{Ncl}(h^{-1}(\mathcal{D}))))) \subseteq \text{Ncl}(\mathcal{D}) = \mathcal{D}$. That is
$\text{Nint}(\text{Ncl}(\text{Nint}(\text{Ncl}(h^{-1}(\mathcal{V} - \mathcal{C}))))) \subseteq h^{-1}(\mathcal{V} - \mathcal{C})$. Or $\text{Nint}(\text{Ncl}(\text{Nint}(\text{Ncl}(h^{-1}(\mathcal{C}))))) \supseteq h^{-1}(\mathcal{C})$. Hence $h^{-1}(\mathcal{C})$ is a NSα-open set in $\mathcal{U}$ and thus $h$ is a NSα-continuous.

**Proposition 3.4.**

i. Every N-continuous map is a Nα-continuous, so it is NSα-continuous, but the opposite is not true in general.
ii. Every Nα-continuous map is a NSα-continuous, but the opposite is not true in general.

*Proof:*

i. Let $h: (\mathcal{U}, \tau_{\mathcal{R}}(\mathcal{M})) \longrightarrow (\mathcal{V}, \sigma_{\mathcal{R}}(\mathcal{N}))$ be a N-continuous map and $\mathcal{D}$ be a N-open set in $\mathcal{V}$. Then $h^{-1}(\mathcal{D})$ is a N-open set in $\mathcal{U}$. Since any N-open set is Nα-open (NSα-open), $h^{-1}(\mathcal{D})$ is a Nα-open (NSα-open) set in $\mathcal{U}$. Thus $h$ is a Nα-continuous (NSα-continuous) map.
ii. Let $h: (\mathcal{U}, \tau_{\mathcal{R}}(\mathcal{M})) \longrightarrow (\mathcal{V}, \sigma_{\mathcal{R}}(\mathcal{N}))$ be a Nα-continuous map and $\mathcal{D}$ be a N-open set in $\mathcal{V}$. Then $h^{-1}(\mathcal{D})$ is a Nα-open set in $\mathcal{U}$. Since any Nα-open set is NSα-open, $h^{-1}(\mathcal{D})$ is a NSα-open set in $\mathcal{U}$. Thus $h$ is a NSα-continuous map.

**Example 3.5.** Let $\mathcal{U} = \{r_1, r_2, r_3, r_4\}$ with $\mathcal{U}/\mathcal{R} = \{\{r_1\}, \{r_4\}, \{r_2, r_3\}\}$ and $\mathcal{M} = \{r_1, r_4\}$. Then $\tau_{\mathcal{R}}(\mathcal{M}) = \{\phi, \{r_1, r_4\}, \mathcal{U}\}$ is a $\text{n.t.s.}$. Let $\mathcal{V} = \{s_1, s_2, s_3, s_4\}$ with $\mathcal{V}/\mathcal{R} = \{\{s_1\}, \{s_3\}, \{s_2, s_4\}\}$ and $\mathcal{N} = \{s_1, s_2\}$. Then $\sigma_{\mathcal{R}}(\mathcal{N}) = \{\phi, \{s_1\}, \{s_2, s_4\}, \{s_1, s_2, s_4\}, \mathcal{V}\}$ is a $\text{n.t.s.}$. Define a map $h: (\mathcal{U}, \tau_{\mathcal{R}}(\mathcal{M})) \longrightarrow (\mathcal{V}, \sigma_{\mathcal{R}}(\mathcal{N}))$ as $h(r_1) = s_2, h(r_2) = s_2, h(r_3) = s_3, h(r_4) = s_4$. Then $h$ is a Nα-continuous but not N-continuous. Also, $h$ is a NSα-continuous but it is not N-continuous.

**Example 3.6.** Let $\mathcal{U} = \{r_1, r_2, r_3, r_4\}$ with $\mathcal{U}/\mathcal{R} = \{\{r_1\}, \{r_3\}, \{r_2, r_4\}\}$ and $\mathcal{M} = \{r_1, r_2\}$. Then $\tau_{\mathcal{R}}(\mathcal{M}) = \{\phi, \{r_1\}, \{r_2, r_4\}, \{r_1, r_2, r_4\}, \mathcal{U}\}$ is a $\text{n.t.s.}$.
Let $\mathcal{V} = \{s_1, s_2, s_3, s_4\}$ with $\mathcal{V}/\mathcal{R} = \{\{s_2\}, \{s_4\}, \{s_1, s_3\}\}$ and $\mathcal{N} = \{s_1, s_2\}$. Then $\sigma_{\mathcal{R}}(\mathcal{N}) = \{\phi, \{s_2\}, \{s_1, s_3\}, \{s_1, s_2, s_3\}, \mathcal{V}\}$ is a $\text{n.t.s.}$. Define a map $h: (\mathcal{U}, \tau_{\mathcal{R}}(\mathcal{M})) \longrightarrow (\mathcal{V}, \sigma_{\mathcal{R}}(\mathcal{N}))$ as $h(r_1) = s_2, h(r_2) = s_1, h(r_3) = s_2,\ h(r_4) = s_3$. It is easily seen that $h$ is a NSα-continuous but it is not Nα-continuous.

**Remark 3.7.** The concepts of N-continuity and Nα*-continuity are independent, for examples.

**Example 3.8.** In example (3.5), the map $h$ is a Nα*-continuous but it is not N-continuous.

**Example 3.9.** Let $\mathcal{U} = \{r_1, r_2, r_3, r_4\}$ with $\mathcal{U}/\mathcal{R} = \{\{r_1\}, \{r_3\}, \{r_2, r_4\}\}$ and $\mathcal{M} = \{r_1, r_2\}$. Then $\tau_{\mathcal{R}}(\mathcal{M}) = \{\phi, \{r_1\}, \{r_2, r_4\}, \{r_1, r_2, r_4\}, \mathcal{U}\}$ is a n.t.s..
Let $\mathcal{V} = \{s_1, s_2, s_3, s_4\}$ with $\mathcal{V}/\mathcal{R} = \{\{s_1\}, \{s_4\}, \{s_2, s_3\}\}$ and $\mathcal{N} = \{s_1, s_4\}$. Then $\sigma_{\mathcal{R}}(\mathcal{N}) = \{\phi, \{s_1, s_4\}, \mathcal{V}\}$ is a n.t.s.. Define a map $h: (\mathcal{U}, \tau_{\mathcal{R}}(\mathcal{M})) \longrightarrow (\mathcal{V}, \sigma_{\mathcal{R}}(\mathcal{N}))$ as $h(r_1) = s_2, h(r_2) = s_1, h(r_3) = s_3, h(r_4) = s_4$. It is easily seen that $h$ is a N-continuous but it is not Nα*-continuous.

**Theorem 3.10.**

i. If a map $h: (\mathcal{U}, \tau_{\mathcal{R}}(\mathcal{M})) \longrightarrow (\mathcal{V}, \sigma_{\mathcal{R}}(\mathcal{N}))$ is N-open, N-continuous and bijective, then $h$ is a Nα*-continuous.
ii. A map $h: (\mathcal{U}, \tau_{\mathcal{R}}(\mathcal{M})) \longrightarrow (\mathcal{V}, \sigma_{\mathcal{R}}(\mathcal{N}))$ is Nα*-continuous iff $h: (\mathcal{U}, \tau_{\mathcal{R}}\alpha(\mathcal{M})) \longrightarrow (\mathcal{V}, \sigma_{\mathcal{R}}\alpha(\mathcal{N}))$ is a N-continuous.

*Proof:*

i. Let $\mathcal{D} \in \sigma_{\mathcal{R}}\alpha(\mathcal{N})$, to prove that $h^{-1}(\mathcal{D}) \in \tau_{\mathcal{R}}\alpha(\mathcal{M})$,
i.e., $h^{-1}(\mathcal{D}) \subseteq \mathrm{Nint}(\mathrm{Ncl}(\mathrm{Nint}(h^{-1}(\mathcal{D}))))$.
Let $a \in h^{-1}(\mathcal{D}) \Longrightarrow h(a) \in \mathcal{D}$. Hence $h(a) \in \mathrm{Nint}(\mathrm{Ncl}(\mathrm{Nint}(\mathcal{D})))$ (since $\mathcal{D} \in \sigma_{\mathcal{R}}\alpha(\mathcal{N})$). Therefore, there exists $\mathcal{H}$ N-open set in $\mathcal{V}$ such that $h(a) \in \mathcal{H} \subseteq \mathrm{Ncl}(\mathrm{Nint}(\mathcal{D}))$. Then $a \in h^{-1}(\mathcal{H}) \subseteq h^{-1}(\mathrm{Ncl}(\mathrm{Nint}(\mathcal{D})))$, but $h^{-1}(\mathrm{Ncl}(\mathrm{Nint}(\mathcal{D}))) \subseteq \mathrm{Ncl}(h^{-1}(\mathrm{Nint}(\mathcal{D})))$ (since $h^{-1}$ is a N-continuous, which is equivalent to $h$ is a N-open and bijective). Then $a \in h^{-1}(\mathcal{H}) \subseteq \mathrm{Ncl}(h^{-1}(\mathrm{Nint}(\mathcal{D})))$. Hence $a \in h^{-1}(\mathcal{H}) \subseteq \mathrm{Ncl}(h^{-1}(\mathrm{Nint}(\mathcal{D}))) \subseteq \mathrm{Ncl}(\mathrm{Nint}(h^{-1}(\mathcal{D})))$ (since $h$ is a N-continuous).
Hence $a \in h^{-1}(\mathcal{H}) \subseteq \mathrm{Ncl}(\mathrm{Nint}(h^{-1}(\mathcal{D})))$, but $h^{-1}(\mathcal{H})$ is a N-open set in $\mathcal{U}$ (since $h$ is a N-continuous). Therefore, $a \in \mathrm{Nint}(\mathrm{Ncl}(\mathrm{Nint}(h^{-1}(\mathcal{D}))))$. Hence $h^{-1}(\mathcal{D}) \subseteq \mathrm{Nint}(\mathrm{Ncl}(\mathrm{Nint}(h^{-1}(\mathcal{D})))) \Longrightarrow h^{-1}(\mathcal{D}) \in \tau_{\mathcal{R}}\alpha(\mathcal{M}) \Longrightarrow h$ is a Nα*-continuous map.
ii. The proof of (ii) is easily.

**Theorem 3.11.** A map $h: (\mathcal{U}, \tau_{\mathcal{R}}(\mathcal{M})) \longrightarrow (\mathcal{V}, \sigma_{\mathcal{R}}(\mathcal{N}))$ is a NSα*-continuous iff $h: (\mathcal{U}, \tau_{\mathcal{R}}\mathrm{S}\alpha(\mathcal{M})) \longrightarrow (\mathcal{V}, \sigma_{\mathcal{R}}\mathrm{S}\alpha(\mathcal{N}))$ is a N-continuous.

*Proof:* Obvious.

**Remark 3.12.** The concepts of N-continuity and NSα*-continuity are independent, for examples:

**Example 3.13.** In example (3.6), the map $h$ is a NSα*-continuous but it is not N-continuous.

**Example 3.14.** Let $\mathcal{U} = \{r_1, r_2, r_3, r_4\}$ with $\mathcal{U}/\mathcal{R} = \{\{r_1\},\{r_4\},\{r_2,r_3\}\}$ and $\mathcal{M} = \{r_1, r_3\}$. Then $\tau_{\mathcal{R}}(\mathcal{M}) = \{\phi, \{r_1\}, \{r_2,r_3\}, \{r_1,r_2,r_3\}, \mathcal{U}\}$ is a n.t.s.. Let $\mathcal{V} = \{s_1, s_2, s_3\}$ with $\mathcal{V}/\mathcal{R} = \{\{s_1\},\{s_2,s_3\}\}$ and $\mathcal{N} = \{s_1, s_3\}$. Then $\sigma_{\mathcal{R}}(\mathcal{N}) = \{\phi, \{s_1\}, \mathcal{V}\}$ is a n.t.s.. Define a map $h: (\mathcal{U}, \tau_{\mathcal{R}}(\mathcal{M})) \longrightarrow (\mathcal{V}, \sigma_{\mathcal{R}}(\mathcal{N}))$ as $h(r_1) = s_1, h(r_2) = s_2$ and $h(r_3) = h(r_4) = s_3$. It is easily seen that $h$ is a N-continuous but it is not NSα*-continuous.

**Remark 3.15.** Every Nα*-continuous map is a Nα-continuous and NSα-continuous but the opposite is not true in general as the following example show:

**Example 3.16.** Let $\mathcal{U} = \{r_1, r_2, r_3, r_4\}$ with $\mathcal{U}/\mathcal{R} = \{\{r_2\},\{r_3\},\{r_1,r_4\}\}$ and $\mathcal{M} = \{r_1, r_3\}$. Then $\tau_{\mathcal{R}}(\mathcal{M}) = \{\phi, \{r_3\}, \{r_1,r_4\}, \{r_1,r_3,r_4\}, \mathcal{U}\}$ is a n.t.s..
Let $\mathcal{V} = \{s_1, s_2, s_3, s_4\}$ with $\mathcal{V}/\mathcal{R} = \{\{s_1\},\{s_2\},\{s_3\},\{s_4\}\}$ and $\mathcal{N} = \{s_1, s_4\}$. Then $\sigma_{\mathcal{R}}(\mathcal{N}) = \{\phi, \{s_1,s_4\}, \mathcal{V}\}$ is a n.t.s.. Define a map $h: (\mathcal{U}, \tau_{\mathcal{R}}(\mathcal{M})) \longrightarrow (\mathcal{V}, \sigma_{\mathcal{R}}(\mathcal{N}))$ as $h(r_1) = s_1, h(r_2) = s_2, h(r_3) = s_3, h(r_4) = s_4$. It is easily seen that $h$ is a Nα-continuous and NSα-continuous but not Nα*-continuous.

**Remark 3.17.** The concepts of Nα*-continuity and NSα*-continuity are independent as the following examples show:

**Example 3.18.** In example (3.16), the map $h$ is a NSα*-continuous but it is not Nα*-continuous.

**Example 3.19.** Let $\mathcal{U} = \{r_1, r_2, r_3, r_4\}$ with $\mathcal{U}/\mathcal{R} = \{\{r_1\},\{r_3\},\{r_2,r_4\}\}$ and $\mathcal{M} = \{r_1, r_2\}$. Then $\tau_{\mathcal{R}}(\mathcal{M}) = \{\phi, \{r_1\}, \{r_2,r_4\}, \{r_1,r_2,r_4\}, \mathcal{U}\}$ is a n.t.s.. Let $\mathcal{V} = \{s_1, s_2, s_3, s_4\}$ with $\mathcal{V}/\mathcal{R} = \{\{s_2\},\{s_4\},\{s_1,s_3\}\}$ and $\mathcal{N} = \{s_1, s_2\}$. Then $\sigma_{\mathcal{R}}(\mathcal{N}) = \{\phi, \{s_2\}, \{s_1,s_3\}, \{s_1,s_2,s_3\}, \mathcal{V}\}$ is a n.t.s.. Define a map $h: (\mathcal{U}, \tau_{\mathcal{R}}(\mathcal{M})) \longrightarrow (\mathcal{V}, \sigma_{\mathcal{R}}(\mathcal{N}))$ as $h(r_1) = h(r_2) = s_1, h(r_3) = s_4$, $h(r_4) = s_3$. It is easily seen that $h$ is a Nα*-continuous but it is not NSα*-continuous.

**Theorem 3.20.** If a map $h: (\mathcal{U}, \tau_{\mathcal{R}}(\mathcal{M})) \longrightarrow (\mathcal{V}, \sigma_{\mathcal{R}}(\mathcal{N}))$ is Nα*-continuous, N-open and bijective, then it is NSα*-continuous.

*Proof:*
Let $h: (\mathcal{U}, \tau_{\mathcal{R}}(\mathcal{M})) \longrightarrow (\mathcal{V}, \sigma_{\mathcal{R}}(\mathcal{N}))$ be a Nα*-continuous, N-open and bijective. Let $\mathcal{D}$ be a NSα-open set in $\mathcal{V}$. Then there exists a Nα-open set say $\mathcal{P}$ such that $\mathcal{P} \subseteq \mathcal{D} \subseteq \mathrm{Ncl}(\mathcal{P})$. Therefore $h^{-1}(\mathcal{P}) \subseteq h^{-1}(\mathcal{D}) \subseteq h^{-1}(\mathrm{Ncl}(\mathcal{P})) \subseteq \mathrm{Ncl}(h^{-1}(\mathcal{P}))$ (since $h$ is a N-open), but $h^{-1}(\mathcal{P}) \in \tau_{\mathcal{R}}\alpha(\mathcal{M})$ (since $h$ is a Nα*-continuous). Hence $h^{-1}(\mathcal{P}) \subseteq h^{-1}(\mathcal{D}) \subseteq \mathrm{Ncl}(h^{-1}(\mathcal{P}))$. Thus, $h^{-1}(\mathcal{D}) \in \tau_{\mathcal{R}}\mathrm{S}\alpha(\mathcal{M})$. Therefore, $h$ is a NSα*-continuous.

**Remark 3.21.** Let $h_1: (\mathcal{U}, \tau_{\mathcal{R}}(\mathcal{M})) \longrightarrow (\mathcal{V}, \sigma_{\mathcal{R}}(\mathcal{N}))$ and $h_2: (\mathcal{V}, \sigma_{\mathcal{R}}(\mathcal{N})) \longrightarrow (\mathcal{W}, \rho_{\mathcal{R}}(\mathcal{O}))$ be two maps, then:

i. If $h_1$ and $h_2$ are Nα-continuous, then $h_2 \circ h_1: (\mathcal{U}, \tau_{\mathcal{R}}(\mathcal{M})) \longrightarrow (\mathcal{W}, \rho_{\mathcal{R}}(\mathcal{O}))$ need not to be a Nα-continuous.
ii. If $h_1$ and $h_2$ are NSα-continuous, then $h_2 \circ h_1: (\mathcal{U}, \tau_{\mathcal{R}}(\mathcal{M})) \longrightarrow (\mathcal{W}, \rho_{\mathcal{R}}(\mathcal{O}))$ need not to be a NSα-continuous.

**Example 3.22.** Let $\mathcal{U}=\{1,2,3,4\}$ with $\mathcal{U}/\mathcal{R}=\{\{2\},\{4\},\{1,3\}\}$ and $\mathcal{M}=\{1,2\}$. Then $\tau_{\mathcal{R}}(\mathcal{M})=\{\phi,\{3\},\{1,3\},\{1,2,3\},\mathcal{U}\}$ is a $\text{n.t.s.}$. The family of all Nα-open (NSα-open) sets of $\mathcal{U}$ is: $\tau_{\mathcal{R}}\alpha(\mathcal{M})=\tau_{\mathcal{R}}\text{S}\alpha(\mathcal{M})=\tau_{\mathcal{R}}(\mathcal{M})\cup\{\{2,3\},\{3,4\},\{1,3,4\},\{2,3,4\}\}$.
Let $\mathcal{V}=\{s_1,s_2,s_3\}$ with $\mathcal{V}/\mathcal{R}=\{\{s_1\},\{s_2\},\{s_3\}\}$ and $\mathcal{N}=\{s_1,s_2\}$. Then $\sigma_{\mathcal{R}}(\mathcal{N})=\{\phi,\{s_3\},\mathcal{V}\}$ is a $\text{n.t.s.}$. The family of all Nα-open (NSα-open) sets of $\mathcal{V}$ is: $\sigma_{\mathcal{R}}\alpha(\mathcal{N})=\sigma_{\mathcal{R}}\text{S}\alpha(\mathcal{N})=\sigma_{\mathcal{R}}(\mathcal{N})\cup\{\{s_1,s_3\},\{s_2,s_3\}\}$.
Define a map $h_1:(\mathcal{U},\tau_{\mathcal{R}}(\mathcal{M}))\longrightarrow(\mathcal{V},\sigma_{\mathcal{R}}(\mathcal{N}))$ as $h_1(1)=h_1(2)=s_1, h_1(3)=h_1(4)=s_2$. Define a map $h_2:(\mathcal{V},\sigma_{\mathcal{R}}(\mathcal{N}))\longrightarrow(\mathcal{U},\tau_{\mathcal{R}}(\mathcal{M}))$ as $h_2(s_1)=h_2(s_3)=3, h_2(s_2)=1$. It is easily seen that $h_1$ and $h_2$ are Nα-continuous (NSα-continuous) maps, but $h_2\circ h_1$: $(\mathcal{U},\tau_{\mathcal{R}}(\mathcal{M}))\longrightarrow(\mathcal{U},\tau_{\mathcal{R}}(\mathcal{M}))$, where $h_2\circ h_1(1)=h_2\circ h_1(2)=3, h_2\circ h_1(3)=h_2\circ h_1(4)=1$, hence $h_2\circ h_1$ is not Nα-continuous (NSα-continuous) map since $\{3\}$ is a N-open set in $\mathcal{U}$, but $(h_2\circ h_1)^{-1}(\{3\})=\{1,2\}$ is not Nα-open (NSα-open) set in $\mathcal{U}$.

**Theorem 3.23.** Let $h_1:(\mathcal{U},\tau_{\mathcal{R}}(\mathcal{M}))\longrightarrow(\mathcal{V},\sigma_{\mathcal{R}}(\mathcal{N}))$ and $h_2:(\mathcal{V},\sigma_{\mathcal{R}}(\mathcal{N}))\longrightarrow(\mathcal{W},\rho_{\mathcal{R}}(\mathcal{O}))$ be two maps, then:

i. If $h_1$ is Nα-continuous and $h_2$ is N-continuous, then $h_2\circ h_1:(\mathcal{U},\tau_{\mathcal{R}}(\mathcal{M}))\longrightarrow(\mathcal{W},\rho_{\mathcal{R}}(\mathcal{O}))$ is a Nα-continuous.
ii. If $h_1$ is Nα*-continuous and $h_2$ is Nα-continuous, then $h_2\circ h_1:(\mathcal{U},\tau_{\mathcal{R}}(\mathcal{M}))\longrightarrow(\mathcal{W},\rho_{\mathcal{R}}(\mathcal{O}))$ is a Nα-continuous.
iii. If $h_1$ and $h_2$ are Nα*-continuous, then $h_2\circ h_1:(\mathcal{U},\tau_{\mathcal{R}}(\mathcal{M}))\longrightarrow(\mathcal{W},\rho_{\mathcal{R}}(\mathcal{O}))$ is a Nα*-continuous.
iv. If $h_1$ and $h_2$ are NSα*-continuous, then $h_2\circ h_1:(\mathcal{U},\tau_{\mathcal{R}}(\mathcal{M}))\longrightarrow(\mathcal{W},\rho_{\mathcal{R}}(\mathcal{O}))$ is a NSα*-continuous.
v. If $h_1$ and $h_2$ are Nα**-continuous, then $h_2\circ h_1:(\mathcal{U},\tau_{\mathcal{R}}(\mathcal{M}))\longrightarrow(\mathcal{W},\rho_{\mathcal{R}}(\mathcal{O}))$ is a Nα**-continuous.
vi. If $h_1$ and $h_2$ are NSα**-continuous, then $h_2\circ h_1:(\mathcal{U},\tau_{\mathcal{R}}(\mathcal{M}))\longrightarrow(\mathcal{W},\rho_{\mathcal{R}}(\mathcal{O}))$ is a NSα**-continuous.
vii. If $h_1$ is Nα**-continuous and $h_2$ is Nα*-continuous, then $h_2\circ h_1:(\mathcal{U},\tau_{\mathcal{R}}(\mathcal{M}))\longrightarrow(\mathcal{W},\rho_{\mathcal{R}}(\mathcal{O}))$ is a Nα**-continuous.
viii. If $h_1$ is Nα**-continuous and $h_2$ is Nα-continuous, then $h_2\circ h_1:(\mathcal{U},\tau_{\mathcal{R}}(\mathcal{M}))\longrightarrow(\mathcal{W},\rho_{\mathcal{R}}(\mathcal{O}))$ is a N-continuous.
ix. If $h_1$ is Nα-continuous and $h_2$ is Nα**-continuous, then $h_2\circ h_1:(\mathcal{U},\tau_{\mathcal{R}}(\mathcal{M}))\longrightarrow(\mathcal{W},\rho_{\mathcal{R}}(\mathcal{O}))$ is a Nα*-continuous.
x. If $h_1$ is N-continuous and $h_2$ is Nα**-continuous, then $h_2\circ h_1:(\mathcal{U},\tau_{\mathcal{R}}(\mathcal{M}))\longrightarrow(\mathcal{W},\rho_{\mathcal{R}}(\mathcal{O}))$ is a Nα**-continuous.

*Proof:*

i. Let $\mathcal{F}$ be a N-open set in $\mathcal{W}$. Since $h_2$ is a N-continuous, ${h_2}^{-1}(\mathcal{F})$ is a N-open set in $\mathcal{V}$. Since $h_1$ is a Nα-continuous, ${h_1}^{-1}({h_2}^{-1}(\mathcal{F}))=(h_2\circ h_1)^{-1}(\mathcal{F})$ is a Nα-open set in $\mathcal{U}$. Thus, $h_2\circ h_1:(\mathcal{U},\tau_{\mathcal{R}}(\mathcal{M}))\longrightarrow(\mathcal{W},\rho_{\mathcal{R}}(\mathcal{O}))$ is a Nα-continuous.
ii. Let $\mathcal{F}$ be a N-open set in $\mathcal{W}$. Since $h_2$ is a Nα-continuous, ${h_2}^{-1}(\mathcal{F})$ is a Nα-open set in $\mathcal{V}$. Since $h_1$ is a Nα*-continuous, ${h_1}^{-1}({h_2}^{-1}(\mathcal{F}))=(h_2\circ h_1)^{-1}(\mathcal{F})$ is a Nα-open set in $\mathcal{U}$. Thus, $h_2\circ h_1:(\mathcal{U},\tau_{\mathcal{R}}(\mathcal{M}))\longrightarrow(\mathcal{W},\rho_{\mathcal{R}}(\mathcal{O}))$ is a Nα-continuous.

iii. Let $\mathcal{F}$ be a Nα-open set in $\mathcal{W}$. Since $h_2$ is a Nα*-continuous, ${h_2}^{-1}(\mathcal{F})$ is a Nα-open set in $\mathcal{V}$. Since $h_1$ is a Nα*-continuous, ${h_1}^{-1}({h_2}^{-1}(\mathcal{F})) = (h_2 \circ h_1)^{-1}(\mathcal{F})$ is a Nα-open set in $\mathcal{U}$. Thus, $h_2 \circ h_1: (\mathcal{U}, \tau_{\mathcal{R}}(\mathcal{M})) \longrightarrow (\mathcal{W}, \rho_{\mathcal{R}}(\mathcal{O}))$ is a Nα*-continuous.

iv. Let $\mathcal{F}$ be a NSα-open set in $\mathcal{W}$. Since $h_2$ is a NSα*-continuous, ${h_2}^{-1}(\mathcal{F})$ is a NSα-open set in $\mathcal{V}$. Since $h_1$ is a NSα*-continuous, ${h_1}^{-1}({h_2}^{-1}(\mathcal{F})) = (h_2 \circ h_1)^{-1}(\mathcal{F})$ is a NSα-open set in $\mathcal{U}$. Thus, $h_2 \circ h_1: (\mathcal{U}, \tau_{\mathcal{R}}(\mathcal{M})) \longrightarrow (\mathcal{W}, \rho_{\mathcal{R}}(\mathcal{O}))$ is a NSα*-continuous.

v. Let $\mathcal{F}$ be a Nα-open set in $\mathcal{W}$. Since $h_2$ is a Nα**-continuous, ${h_2}^{-1}(\mathcal{F})$ is a N-open set in $\mathcal{V}$. Since any N-open set is a Nα-open, ${h_2}^{-1}(\mathcal{F})$ is a Nα-open set in $\mathcal{V}$. Since $h_1$ is a Nα**-continuous, ${h_1}^{-1}({h_2}^{-1}(\mathcal{F})) = (h_2 \circ h_1)^{-1}(\mathcal{F})$ is a N-open set in $\mathcal{U}$. Thus, $h_2 \circ h_1: (\mathcal{U}, \tau_{\mathcal{R}}(\mathcal{M})) \longrightarrow (\mathcal{W}, \rho_{\mathcal{R}}(\mathcal{O}))$ is a Nα**-continuous. The proof is obvious for others.

**Remark 3.24.** The following diagram explains the relationship between different classes of weakly nano continuous maps:

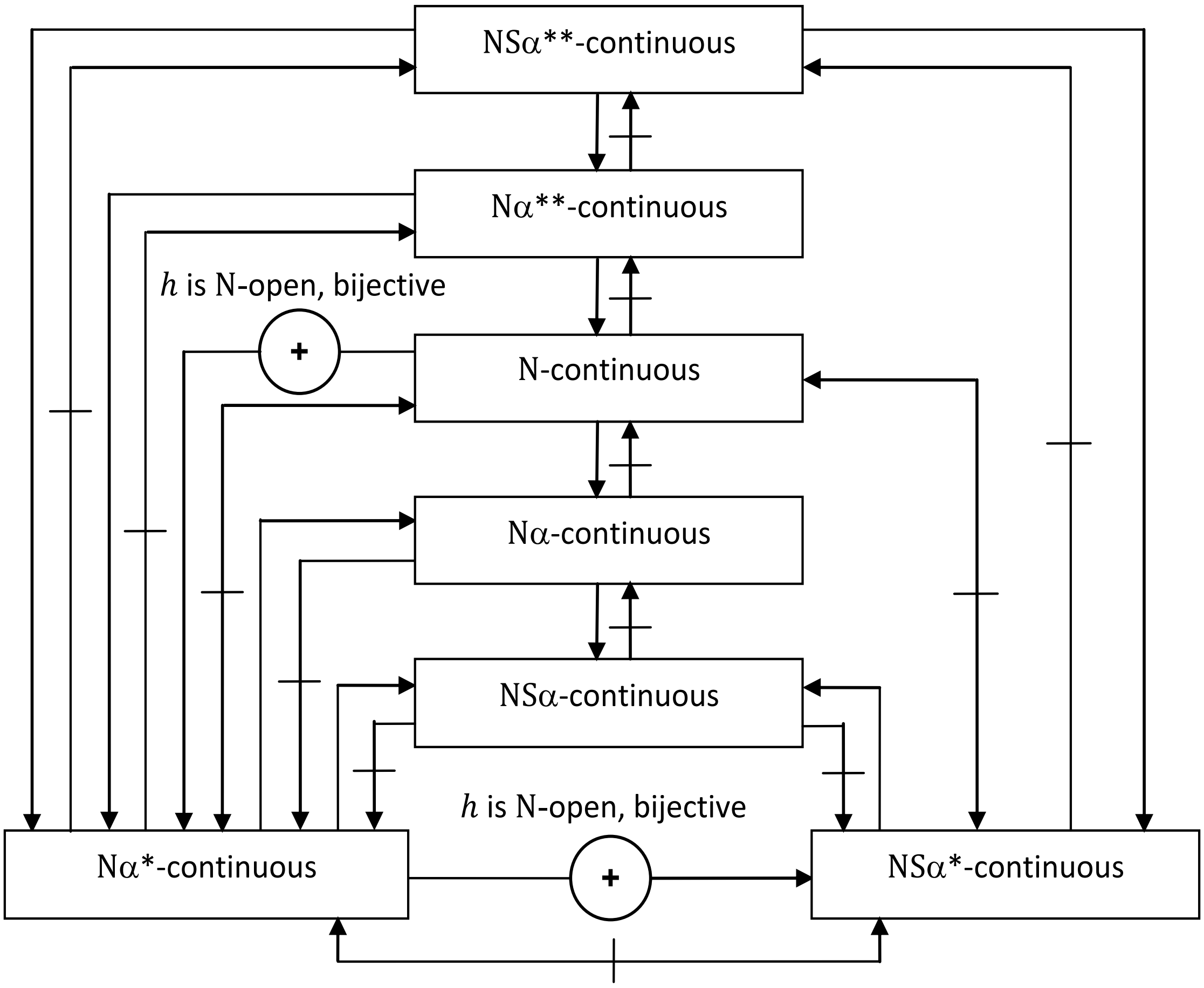

**Diagram (3.1)**

## 4. CONCLUSION

We shall use the concepts of Nα-open and NSα-open sets to define some new types of weakly nano continuity such as; Nα-continuous, Nα*-continuous, Nα**-continuous, NSα-continuous, NSα*-continuous and NSα**-continuous maps. The Nα-open and NSα-open sets can be used to derive some new types of weakly nano open maps, nano compactness, and nano connectedness.